\documentclass[12pt]{article}
\usepackage{amsmath,amsfonts,amssymb}
\usepackage{latexsym}

\def\1{\mbox{\bf 1}}

\def\a{\alpha}
\def\b{\beta}
\def\c{\cdot}
\def\cs{\cdots}
\def\r{\rightarrow}
\def\R{\mathbb{R}}
\def\P{\mathbb{P}}
\def\E{\mathbb{E}}

\def\z{\zeta}
\def\t{\tau}
\def\o{\omega}
\def\iy{\infty}
\def\qed{\hfill\vrule height5pt width5pt depth0pt}

\newcommand{\bea}{\begin{eqnarray}}
\newcommand{\eea}{\end{eqnarray}}
\newcommand{\Bea}{\begin{eqnarray*}}
\newcommand{\Eea}{\end{eqnarray*}}
\newcommand{\proof}{\noindent {\bf Proof:\ }}

\newtheorem{Definition}{Definition}[section]
\newtheorem{Theorem}[Definition]{Theorem}
\newtheorem{Lemma}[Definition]{Lemma}
\newtheorem{Proposition}[Definition]{Proposition}
\newtheorem{Corollary}[Definition]{Corollary}
\newtheorem{Remark}[Definition]{Remark}
\newtheorem{Example}[Definition]{Example}

\begin{document}

\title{ Duality for multidimensional ruin problem }

\author{
S. Ramasubramanian\thanks{Theoretical Statistics and Mathematics Unit; Indian Statistical Institute; 8th mile, Mysore Road; Bangalore - 560 059; India ({\tt ram@isibang.ac.in}).}
}
\date{August 09, 2014}
\maketitle
\begin{abstract}
We consider a $d-$dimensional insurance network, with initial capital $a\in\R^d_+,$ operating under a risk diversifying treaty; this is described in terms of a regulated random walk $\{Z^{(a)}_n\}$ via Skorokhod problem in $\R^d_+$ with reflection matrix $R;$ $\{Y^{(a)}_n\}$ denotes the corresponding pushing process. Ruin (in a strong sense) of $\{Z^{(a)}_n\}$ is defined as the marginal deficit of each company being positive (and hence zero surplus) at some time $n.$ A dual storage network is introduced through time reversal at sample path level over finite time horizon; the stochastic analogue is again a regulated random walk $\{W_n\}$ in $\R^d_+$ starting at $0.$ 
It is shown that ruin for $\{Z^{(a)}_n\}$ corresponds to $\{W_n\}$ hitting open upper orthant determined by $R^{-1}a$ before hitting the boundary of $\R^d_+,$ even at the sample path level. 
Under natural hypotheses, we show that $\P($ ruin of $\{Z^{(a)}_n\}$ in finite time) $=\lim_{n\r\iy}\P(W_n\gg R^{-1}a:~n<$ boundary hitting time of storage process) $=\lim_{n\r\iy}\P(Y^{(0)}_n \gg R^{-1}a:\Delta Y^{(0)}_n\gg 0).$ A notion of $d-$dimensional ladder height distribution is defined, and a Pollaczek-Khinchine formula derived; an expression for the ladder height distribution is presented. Our method is applicable to ruin problem for a continuous time $d-$dimensional Cramer-Lundberg type network, where the companies act independently in the absence of treaty.

\end{abstract}

\section{ Introduction}
\setcounter{equation}{0}

Connection between ruin probabilities of actuarial risk theory and asymptotic behaviour of storage processes in queuing theory is well known in the one dimensional context for more than fifty years; for example, see \cite{Pr1961,Sg1976}. Such a connection has been inspired by the so called duality results for random
walks in $\R;$ see \cite{Sp1956}, and especially Chapter XII of \cite{F1969} and \cite{Bi2001}.
 A comprehensive exposition is given in \cite{AA2010}.

In recent years there has been considerable interest in multidimensional insurance models, including the associated ruin problems; see \cite{APP2008,BG2008,Co1998,CYZ2003} for example. Notions of ruin in these relate to the vector current surplus process hitting a preassigned open set; the latter is generally taken to be the complement of the nonnegative orthant (corresponding to at least one insurance company in the network getting ruined), or the negative orthant (corresponding to all the companies getting ruined at the same time), or a preassigned half space (corresponding to the current total surplus of all companies going below a certain level).

As pointed out by Buhlman, in spite of the unfortunate terminology, the term `ruin of a company' does not imply that the concerned company is crashing out of business, but only highlights a ``need for additional capital"; see p.133 of \cite{Bu1970}. It is also referred to as ``capital injection by the shareholders of the company" in \cite{DW2004}. Well known optimality properties of one and higher dimensional Skorokhod problem (see \cite{Ha1985,Re1984,CM1991,Ra2000}) suggest an optimal way of going about it.

A few years back we had proposed a multidimensional insurance model in terms of Skorokhod problem (SP, for short) in an orthant, describing the joint dynamics of $d$ insurance companies operating under a risk diversifying treaty. According to the treaty, when a company
in the network needs an amount to prevent its surplus from getting
wiped out, the required capital injection is obtained from other
companies in the network, as well as from the shareholders in pre-agreed proportions; and the optimal way to go about
is provided by the SP; see the discussion just after Theorem~\ref{Th:splcp} below. The reflection matrix will not be diagonal in general, that is, we need to consider oblique reflection. It has been argued in \cite{Ra2006,Ra2011} that it results in a reasonable model. So the regulated/ reflected part of the solution to the SP gives the optimal (vector) current surplus, and the pushing part of the solution gives the optimal (vector) capital injection (for averting ruin), while operating under the risk diversifying treaty.

In this paper we consider the ruin problem for multidimensional insurance models that are described in terms of regulated random walks in a $d-$dimensional orthant. The reflection matrix is taken to be constant matrix. Clearly, the vector $0$ has a special status, and this leads to canonical notion(s) of ruin of the network; see also \cite{Ra2012}. We give $3$ natural, but closely related, notions of ruin: ss-ruin corresponding to  each company needing positive capital injection, that is, each company having nonzero marginal deficit (and hence zero surplus), at time $n,$ s-ruin corresponding to each company having zero surplus with at least one having nonzero marginal deficit as well at $n,$ and ruin corresponding to each company having zero surplus at $n,$ for some $n.$ Note that all are  connected to the surplus process hitting the state $0$ in finite time. Under minimal conditions, these three notions coincide with probability $1;$ besides, in the one-dimensional case, these coincide with the classical notion of ruin.

In our setup, the SP for the sequence of partial sums can be built out of a sequence of Linear Complementarity problems for a sequence of vectors, of course, corresponding to the same reflection matrix $R.$ When the (discrete) time horizon is finite, through time reversal, we are led naturally to a dual discrete time regulated random walk in the $d-$dimensional orthant, which is referred to as a storage network. This storage network admits a reasonable interpretation. In finite time horizon, ruin of the insurance network is characterized in terms of the dual storage network crossing a certain threshold, at the sample path level; the matrix $R^{-1}$ plays a major role. For considering the stochastic setup in infinite time horizon, we introduce various hypotheses, including the coordinatewise net profit condition. It is shown that the ss-ruin probability of the insurance network can be expressed as the probability of the storage network exceeding a certain threshold (given in terms of the initial capital) \emph{before} hitting the boundary of the orthant. Moreover, it turns out that the asymptotic behaviour of the dual storage network before hitting the boundary, and the asymptotic behaviour of the pushing process (when it is strictly increasing) associated with the insurance network (with initial capital $0$) are closely related. 
We also introduce an appropriate notion of $d-$dimensional ladder height distribution, and obtain a Pollaczek-Khinchine formula for ss-ruin probability; we are able to express the ladder height distribution in terms of the given data.

We now indicate a class of examples covered by our analysis. Suppose that, in the absence of the risk diversifying treaty, the joint dynamics of the companies is a continuous time $d-$dimensional renewal risk process given by \eqref{Hcts} in Example~\ref{Ex:renew}. The scalar i.i.d. interarrival times, the random mechanism governing which among the companies would take the claim at an arrival time, and the i.i.d. $d-$dimensional claim size vectors form independent families of random variables. To study the ruin problem in this case it is enough to consider the process at claim arrival times; and the process observed only at claim arrival times constitutes a random walk in $\R^d.$ So our method is applicable to study the ruin problem for such processes. An important special case is that of a Cramer-Lundberg type network; in the absence of the treaty, the joint dynamics is that of $d$ independent one dimensional Cramer-Lundberg processes.

To our knowledge \cite{BS1999} seems to be the only other paper to have considered duality and multidimensional risk models. However, the approach and emphasis seem to be quite different from ours. For example, in \cite{BS1999}, the  queueing process and the dual risk process may be based on spaces whose dimensions widely differ, with the latter being set-valued in general; also only normal reflection has been considered.

We now briefly outline the organisation of the paper. Section 2 deals with the deterministic setup, while Section 3 concerns the stochastic setup. In Section 2, we introduce insurance networks described in terms of regulated random walks in an orthant, and the notion of ruin for such networks. The  dual discrete time storage network for an insurance network is then presented, over a finite time horizon. This section concludes with sample-path characterization of ruin in terms of dual storage network. Stochastic analogues are considered in Section 3 along with appropriate hypotheses. Duality results, in the sense of equality in distribution, are derived. A Pollaczek-Khinchine formula for ss-ruin probability is obtained, and the ladder height distribution is identified using duality. A detailed discussion concerning ruin problem for renewal risk type network is also given. 

We now conclude Section 1 with the list of all hypotheses needed in the sequel.

\subsection{ Hypotheses }

Notation: We shall denote by SP and LCP, respectively, the Skorokhod problem and the linear complementarity problem. For $x\in\R^d,$ $(x)_i$ denotes the $i$-th component of $x.$   For $x,y\in\R^d,$ we shall write: $x\geq y$ if $(x)_i\geq (y)_i$ for all $1\leq i\leq d;$ $x>y$ if $x\geq y$ with $(x)_i>(y)_i$ at least for some $i;$ $x\gg y$ if $(x)_i>(y)_i$ for all $1\leq i\leq d.$ Also for vectors $x,y,$ $x\ll y$ is the same as $y\gg x;$ similar comments apply to $ x\leq y, x < y.$

$\overline{G}\triangleq\R^d_+$  denotes the $d$-dimensional nonnegative orthant, and $G$ denotes its interior $\{x\in\R^d: x\gg 0\}.$ All random variables and processes are defined on a probability space $(\Omega,\mathcal{F},\P);$ $\E$ denotes expectation w.r.t. $\P.$

Vectors will be denoted by lower case alphabets, while random variables by capital letters.

\begin{description}

\item{\textbf {(H1)}} $R=((R_{ij}))=I-P^t$ is a $d\times d$ real matrix such that $P_{ii}=0,~P_{ij}\geq 0,~i\neq j,$ for all $1\leq i,j\leq d,$  and spectral radius of $P$ is strictly less than $1.$  Here $R$ denotes reflection matrix.

\item{\textbf {(H2)}} There exists $k\in \{1,2,\cs,d\}$ such that $(R^{-1})_{ik}>0$ for all $1\leq i\leq d;$ that is, at least one column vector of $R^{-1}$ has strictly positive entries.

\item{\textbf{(H3)}} $A_i,i=1,2,\ldots $ denote one dimensional i.i.d. random variables such that $A_i > 0;$ these are (scalar) interarrival times.

\item{\textbf{(H4)}} $X_{\ell},\ell =1,2,\ldots $ are i.i.d. $\R^d_+$-valued random variables; these are vector claim sizes.

\item{\textbf {(H5)}} $\{A_i:i\geq 1\},\{X_{\ell}:\ell \geq 1\}$ are independent families of random variables.

\item{\textbf {(H6)}} For each $\ell = 1,2,\cs $ and $i=1,2,\cs,d,$ $\P((X_{\ell})_i >x) >0,~~\forall x\geq 0;$ that is, marginal claim sizes have unbounded support.

\item{\textbf {(H7)}} For each $\ell = 1,2,\cs $ and $i=1,2,\cs,d,$ $\P((X_{\ell})_i = x) = 0,~~\forall x > 0;$ that is, $(X_{\ell})_i$ has no atoms in $(0,\iy);$ however, there can be an atom at $0.$

\item{\textbf {(H8)}} $c=((c)_1,\cs,(c)_d) \gg 0$ with $(c)_i$ denoting constant premium rates. $A_1,(X_1)_i,1\leq i\leq d$ have finite expectations, and $\E[(c)_iA_1 - (X_1)_i] > 0,~~1\leq i\leq d;$ this is coordinatewise net profit condition.

\end{description}

Note that (H1),(H2) concern only the reflection matrix $R,$ and involve no probabilistic assumptions. Our  analysis on deterministic set up in Section 2 will involve only (H1).

\begin{Remark} {\rm  (i) By the spectral radius condition in (H1) note that \bea\label{Rinv} R^{-1} &=& I+P^t+(P^t )^2+(P^t)^3+\cs \eea is a matrix with nonnegative entries, with diagonal entries $\geq 1.$

(ii) In the context of insurance models, in addition to (H1), it is natural to assume that $\sum_{j\neq i} P_{ij} \leq 1$ for all $i,$ that is $P$ is a substochastic matrix.

(iii) Note that (H2) holds if $P$ is irreducible; see \cite{Se1981}. It also holds in the feedforward case.
}
\end{Remark}

\section{ Deterministic setup }
\setcounter{equation}{0}

In this section we introduce the deterministic analogues of insurance and storage networks described in terms of regulated random walks in an orthant. We establish duality results in a finite discrete time horizon at sample path level.

\subsection{ SP and LCP }

We now describe Skorokhod problem (SP, for short) in an orthant for partial sums in the deterministic set up; this basically involves solving a sequence of linear complementarity problems (LCP). Required references on SP will be given at appropriate places, while \cite{CPS1992} is an encyclopaedic work on LCP; \cite{Ma1989} gives an exposition on the connection between SP and LCP.

Let $R$ be a reflection matrix satisfying (H1). Let $a=((a)_1,\cs,(a)_d)\in \overline{G}.$ Let $\{u_n,~n\geq 1\}$ denote a sequence in $\R^d.$ A pair $\{y^{(a)}_n,~n\geq 0\},\{z^{(a)}_n,~n\geq 0\}$ of sequences in $\R^d$ is said to be a solution to the \emph{deterministic Skorokhod problem} $SP(\{a+\sum u_n\},R)$ if the following hold:

\begin{description}

\item{\bf (s0)} $y^{(a)}_0=0,~~z^{(a)}_0=a.$

\item{\bf (s1)} For $1\leq i\leq d$ Skorokhod equation holds, that is, \bea\label{seqi} (z^{(a)}_n)_i &=& (a)_i + \sum_{\ell =1}^n (u_{\ell})_i + (y^{(a)}_n)_i  +  \sum_{j\neq i} R_{ij}(y^{(a)}_n)_j,~~n\geq 1; \eea or equivalently in vector notation \bea\label{seq} z^{(a)}_n &=& a + \sum_{\ell =1}^n u_{\ell} + Ry^{(a)}_n \nonumber \\ & = & z^{(a)}_{n-1} + u_n + R\Delta y^{(a)}_n,~~n\geq 1, \eea where $\Delta y^{(a)}_n = y^{(a)}_n - y^{(a)}_{n-1}.$

\item{\bf (s2)} $z^{(a)}_n\in \overline{G}$ for all $n\geq 1.$

\item{\bf (s3)} $y^{(a)}_n \geq y^{(a)}_{n-1},~n\geq 1$ as vectors; moreover $(y^{(a)}_{\c})_i$ can increase only when $(z^{(a)}_{\c})_i=0,$ that is  \bea\label{mincondn} \langle z^{(a)}_n,\Delta y^{(a)}_n \rangle &=& 0,~~n\geq 1. \eea

\end{description}
Note that (s2) is a constraint, while \eqref{mincondn} in (s3) is a minimality condition. We refer to  $\{y^{(a)}_n\},\{z^{(a)}_n\}$ respectively as the \emph{pushing part, regulated/ reflected part} of the solution to $SP(\{a+\sum u_n\},R).$

To describe the \emph{linear complementarity problem}, let $\eta\in\R^d$ and $R$ as above. We say a pair $\xi,\z\in\R^d$ is a solution to $LCP(\eta,R)$ if $\z=\eta + R\xi,$ $\xi\geq 0,~\z\geq 0$ as vectors, and $\langle \xi,\z\rangle =0.$ We denote $\xi=\Phi(\eta,R),~\z=\Psi(\eta,R)$ and call them respectively the pushing part, regulated part of the solution.

A compilation of required results on deterministic (or equivalently sample path) SP for partial sums is given below. For details and proofs see \cite{HR1981,Re1984,Ma1989,CM1991,KW1996,R2006,Ra2000} and references therein.

\begin{Theorem}\label{Th:splcp} Let $R$ satisfy (H1); let $u_n\in\R^d,~n\geq 1.$ Then the following hold.

(i) There is a unique solution pair to $LCP(\eta,R)$ for any $\eta\in\R^d.$

(ii) There is a unique solution pair  $\{y^{(a)}_n\},\{z^{(a)}_n\}$ to $SP(\{a+\sum u_n\},R)$ for any $a\in\overline{G}.$

(iii) $\{y^{(a)}_n\},\{z^{(a)}_n\}$ is the solution pair to $SP(\{a+\sum u_n\},R)$ if and only if for $n=1,2,\cs $ $(\Delta y^{(a)}_n,z_n)$ is the solution to $LCP(z^{(a)}_{n-1} + u_n,R).$

(iv) If $a,b\in\overline{G}$ with $a\leq b,$ then for $n\geq 1,$ \bea\label{yineq} \Delta y^{(a)}_n & \geq & \Delta y^{(b)}_n, \\ \label{zineq} z^{(a)}_n & \leq & z^{(b)}_n, \\ \label{ybd} 0\leq y^{(a)}_n-y^{(b)}_n & \leq & R^{-1}(b-a). \eea

(v) For $n\geq 1,$ put $h^{(a)}_n=((h^{(a)}_n)_1,\cs,((h^{(a)}_n)_d)$ where \Bea (h^{(a)}_n)_i & = & \sup_{k\leq n} \max \{0,-((a)_i + \sum_{\ell =1}^k (u_{\ell})_i)\}.\Eea Then \bea\label{apri} y^{(a)}_n & \leq & R^{-1} h^{(a)}_n,~~n\geq 1. \eea

\end{Theorem}

See Theorem 6 of \cite{KW1996}, and Proposition 3.2 and Theorem 4.1 of \cite{Ra2000} for proofs of \eqref{yineq} - \eqref{apri}.

The framework above can be given the following interpretation in the context of insurance models. We consider $d$ insurance companies operating under a risk diversification treaty specified in terms of $R,$ with $P$ being a substochastic matrix. Claims are assumed to arrive at regular intervals $k=1,2,\cs,$ and are settled instantaneously. According to the treaty, if Company $i$ requires an amount $(\Delta \theta_k)_i$ at time $k$ to avert ruin, then for $j\neq i,$ Company $j$ gives $P_{ij}(\Delta \theta_k)_i = |R_{ji}|(\Delta \theta_k)_i$ from its surplus; any shortfall has to be provided by the shareholders of Company $i$ as capital injection; for $j\neq i,$ if Company $j$ is not able give from its surplus, then Company $j$ may also have to get capital injection. So the surplus of any company is required to be nonnegative. The spectral radius condition in (H1) means that the network is `open', in the sense that capital injection for the entire network is also possible; this makes the Skorokhod problem well posed. With each company striving to minimize its liability, Skorokhod problem provides the optimal way of operating under the treaty. So under optimality, a company can invoke the treaty only when it is in the red, and the amount it gets from all sources is just enough to keep it afloat. See \cite{Ra2006,Ra2011}. In view of the above, note that

$(a)_i=$ initial capital of Company $i;$

$(u_k)_i=$ (premium income for Company $i$ during $(k-1,k])$ minus (claim amount for Company $i$ due to $k-$th claim);

$(z^{(a)}_k)_i=$ current surplus of Company $i$ at time $k,$ under optimality;

$(y^{(a)}_k)_i=$ cumulative amount obtained by Company $i$ from all sources till time $k$ for the purpose of averting ruin, under optimality; so $\Delta y^{(a)}_k)_i=$ marginal deficit of Company $i$ at time $k,$ under optimality.

Thus the regulated/ reflected part of the solution to $SP(\{a+\sum u_n\},R)$ gives the optimal joint dynamics of $d$ companies operating under the treaty. Many notions introduced in the sequel are related to insurance models.

\subsection{ 3 notions of ruin }

Clearly the state $0$ has a special place in our set up. In \cite{Ra2012} we had defined ruin as the event that the regulated process hitting the origin; this definition works well when dealing with continuous random variables. However this definition may not be strong enough if the process can hit $0$ without any sector being in deficit. Therefore we now define three natural closely related notions of ruin; as we shall see later, these notions coincide under certain probabilistic assumptions.

Let $a,R,\Delta y^{(a)}_n,z^{(a)}_n$ be as above. We say \emph{ruin} occurs for $\{z^{(a)}_n\}$ if $z^{(a)}_k=0$ for some $k\geq 1;$ similarly \emph{s-ruin} occurs for $\{z^{(a)}_n\}$ if $z^{(a)}_k=0,~\Delta y^{(a)}_k >0$ for some $k\geq 1;$ and  \emph{ss-ruin} occurs for $\{z^{(a)}_n\}$ if $\Delta y^{(a)}_k \gg 0$ for some $k\geq 1.$

Because of the minimality condition \eqref{mincondn}, note that ss-ruin implies s-ruin which in turn implies ruin. Note that ss-ruin denotes each company having nonzero marginal deficit (and hence zero surplus) at time $n,$ while s-ruin corresponds to each company having zero surplus with at least one having nonzero marginal deficit as well at time $n,$ and ruin means that each company having zero surplus at time $n,$ for some $n.$

Let $n\geq 1$ be fixed. Using \eqref{mincondn} for $n,$ the Skorokhod equation \eqref{seq} successively for $k\leq n,$ and uniqueness of the solution to $LCP(z^{(a)}_{n-1}+u_n,R),$ we get \bea\label{ssruinequi} \Delta y^{(a)}_n \gg 0 & \Leftrightarrow & -R^{-1}u_n \gg R^{-1}z^{(a)}_{n-1} \nonumber\\ & \Leftrightarrow  & -R^{-1}u_n -R^{-1}u_{n-1} \gg R^{-1}z^{(a)}_{n-2} + \Delta y^{(a)}_{n-1} \nonumber \\ & \Leftrightarrow & \sum_{\ell = n-k}^n (-R^{-1}u_{\ell}) \gg R^{-1}z^{(a)}_{n-(k+1)} + [y^{(a)}_{n-1} - y^{(a)}_{n-(k+1)}] \nonumber \\ & \Leftrightarrow & \sum_{\ell=2}^n (-R^{-1}u_{\ell}) \gg R^{-1}z^{(a)}_1 + [y^{(a)}_{n-1}-y^{(a)}_1] \nonumber \\ & \Leftrightarrow & -R^{-1}a + \sum_{\ell=1}^n (-R^{-1} u_{\ell}) \gg y^{(a)}_{n-1} \eea (Note that the difference between the two sides of the last inequality in the string \eqref{ssruinequi} is $\Phi(z^{(a)}_{n-1}+u_n,R),$ by uniqueness of solution to LCP.) Thus we have
\begin{Proposition}\label{Prop:3ruin} Let $R$ satisfy (H1); let $n\geq 1$ be fixed. Then the following hold.
\bea\label{ssruin} \Delta y^{(a)}_n \gg 0 & \Leftrightarrow & -R^{-1}u_n \gg R^{-1}z^{(a)}_{n-1} \nonumber\\ & \Leftrightarrow  & -R^{-1}a + \sum_{\ell=1}^n (-R^{-1} u_{\ell}) \gg y^{(a)}_{n-1}. \\ \label{sruin}
z^{(a)}_n=0,~\Delta y^{(a)}_n > 0 & \Leftrightarrow & -R^{-1}u_n > R^{-1}z^{(a)}_{n-1} \nonumber\\ & \Leftrightarrow  & -R^{-1}a + \sum_{\ell=1}^n (-R^{-1} u_{\ell}) > y^{(a)}_{n-1}. \\ \label{ruin} z^{(a)}_n=0 & \Leftrightarrow & -R^{-1}u_n \geq R^{-1}z^{(a)}_{n-1} \nonumber\\ & \Leftrightarrow  & -R^{-1}a + \sum_{\ell=1}^n (-R^{-1} u_{\ell}) \geq y^{(a)}_{n-1}. \eea Moreover, in all three cases, if left side holds then \bea\label{yan} y^{(a)}_n & = & -R^{-1}a + \sum_{\ell=1}^n (-R^{-1} u_{\ell}). \eea
\end{Proposition}
\proof \eqref{ssruin} has been established prior to the statement of the proposition. Replacing $\gg$ by $ >, \geq,$ respectively, in the proof of \eqref{ssruin}, one can prove \eqref{sruin}, \eqref{ruin}. Now \eqref{yan} is an easy consequence.  \qed

\subsection{Storage network in finite discrete time horizon}

We begin with an elementary observation.
\begin{Proposition}\label{Prop:lcphat} Let $R$ satisfy (H1). Then for any $\eta\in\R^d,$  \bea\label{lcpR-1}  & & \chi=\Phi(\eta,R),~\varpi=\Psi(\eta,R) \nonumber \\ & \Leftrightarrow & \varpi=\Phi(-R^{-1}\eta,R^{-1}),~\chi=\Psi(-R^{-1}\eta,R^{-1}). \eea In particular, $LCP(\theta,R^{-1})$ has a unique solution pair for any $\theta\in\R^d.$
\end{Proposition}
\proof Clearly \Bea \varpi=\eta+R\chi & \Leftrightarrow & \chi = -R^{-1}\eta + R^{-1}\varpi. \Eea Hence, whenever $\chi\geq 0,\varpi\geq 0, \langle \chi,\varpi \rangle =0 $ hold, \eqref{lcpR-1} would also hold. As $R$ is invertible, uniqueness of LCP corresponding to $R^{-1}$ follows from that of LCP corresponding to $R.$  \qed

We now consider Skorokhod problem for a collection of partial sums related to the earlier one through time reversal. To describe the sample path (or equivalently the deterministic) set up, we need to look at a finite discrete time horizon.

Assume that $R$ satisfies (H1). Let $n\geq 1$ be fixed. Let $u_k\in\R^d,~1\leq k\leq n.$ Set $\hat{u}_1=-R^{-1}u_n,~\hat{u}_2=-R^{-1}u_{n-1},\cs,~\hat{u}_n=-R^{-1}u_1;$ so $\hat{u}_k=-R^{-1}u_{n+1-k},~1\leq k\leq n.$ Put $w_0=0,~v_0=0.$ For $1\leq k\leq n,$ let $\Delta v_k,w_k$ be the unique solution guaranteed by Proposition~\ref{Prop:lcphat} to $LCP(w_{k-1}+\hat{u}_k,R^{-1}).$ So $\Delta v_k\geq 0,~w_k\geq 0,~\langle w_k,\Delta v_k\rangle =0,$ and \bea\label{seqhat} w_k &=& \sum_{\ell =1}^k \hat{u}_{\ell} + R^{-1}v_k \nonumber \\ & = & w_{k-1} + \hat{u}_k + R^{-1}\Delta v_k,~ 1\leq k\leq n \eea where $v_k=v_0+\sum_{\ell=1}^k\Delta v_{\ell}.$ That is, in the spirit of Theorem~\ref{Th:splcp}, $v_k,w_k,~1\leq k\leq n,$ solve the Skorokhod problem $SP(\{\sum_{\ell=1}^k\hat{u}_{\ell},1\leq k\leq n\},R^{-1}).$ Note that uniqueness of the solution to $SP(\{\sum_{\ell=1}^k\hat{u}_{\ell},1\leq k\leq n\},R^{-1})$ follows from   Proposition~\ref{Prop:lcphat}; see also \cite{Ma1989}. We refer to this set up as a deterministic storage network in finite discrete time horizon; here, $v_k,w_k$ are, respectively, pushing and regulated parts of the storage network.

We now give an interpretation of the storage network. Suppose there are $d$ storage depots of infinite capacity; let the initial stock be $0$ at each depot. While demands might be continuously made, fresh stocks and reinforcements arrive only at the end of periods $k=1,2\cs;$ readings only at the end of the periods are available. Following assumptions are made.

(a) All demand at a depot during a certain period is met at the end of the same period, if necessary by bringing in reinforcement.

(b) A need for reinforcement at Depot $i$ at the end of period $k,$  indicates that available stock at the end of period $k,$ including the arrival (and possible inflow as given in this paragraph later) at the end of period $k,$ has not been sufficient to fulfil the demand. This can trigger increased demand at Depots $j\neq i$ during subsequent periods. So reinforcements are sent to Depots $j\neq i$ as well at the end of the same period $k;$ such an inflow at Depot $j$ can also be used to take care of possible unfulfilled demand at that depot at the end of period $k.$ (Such a mechanism may be motivated by a desire to avoid wider customer dissatisfaction in a cooperative setting, or as an attractive business opportunity in a competitive setting.)

(c) Reinforcement supplied to Depot $i$ at the end of period $k$ (due to unfulfilled demand at that depot) is just enough to fulfil the shortfall at the end period $k;$ that is, reinforcement is `minimal'.

The above interpretation leads to the following meanings.

$(\hat{u}_k)_i=$ (amount of fresh supply arriving at Depot $i$ at the end of period $k$) minus (demand at Depot $i$ during the period $(k-1,k]$);

$(w_k)_i=$  current stock at Depot $i$ at the end of period $k,$ after taking into account all reinforcements to Depot $i$ till the end of period $k;$ so $(w_k)_i\geq 0$ for all $i,k;$

$(R^{-1})_{ii}(\Delta v_k)_i=$ amount of reinforcement sent to Depot $i$ at the end of period $k,$ due to unfulfilled demand after taking into account existing stock, fresh supply and inflow to Depot $i$ due to shortfall at other depots at the end period $k;$

$(R^{-1})_{ij}(\Delta v_k)_j=[(R^{-1})_{ij}/(R^{-1})_{jj}](R^{-1})_{jj}(\Delta v_k)_j=$ amount of reinforcement (inflow) sent to Depot $i$ due to shortfall at Depot $j,$ for $j\neq i,$ at the end of period $k.$

Therefore note that \Bea (\Delta v_k)_i > 0 &\Leftrightarrow & (w_{k-1})_i+(\hat{u}_k)_i +\sum_{j\neq i}(R^{-1})_{ij}(\Delta v_k)_j < 0  \\ & \Leftrightarrow & -(\hat{u}_k)_i > (w_{k-1})_i +\sum_{j\neq i}(R^{-1})_{ij}(\Delta v_k)_j. \Eea
In such a case $(w_k)_i=0,$ that is, \bea\label{hat-v} (R^{-1})_{ii}(\Delta v_k)_i & = & -[(w_{k-1})_i+(\hat{u}_k)_i +\sum_{j\neq i}(R^{-1})_{ij}(\Delta v_k)_j]. \eea

{\bf Note:} The storage network described above might be suitable when the depots are viewed upon as different banks in a small geographical region. Reinforcement at one bank can result in (defensive) inflow at other banks; of course, it is assumed that the exact quantum of reinforcement at one bank is known (or made known) to other banks without delay. The set up can also be looked upon as different branches of the same bank, with reinforcements coming only from a central node (which is not considered part of the network).

\subsection{ A connection }

Let $n\geq 1$ be fixed; let $u_k\in\R^d,~1\leq k\leq n,~a\in\overline{G},$ and $R$ be a matrix as before.  We consider $\{y^{(a)}_k,z^{(a)}_k:1\leq k\leq n\}$ and $\{v_k,w_k:1\leq k\leq n\}$ defined earlier.

Define \bea\label{hittimebd} \sigma_{bd} &=& \inf\{k\geq 1: w_k\in\partial G \}, \\
 \label{exitss} \vartheta_{R^{-1}a} &=& \inf\{k\geq 1: w_k \gg R^{-1}a \}; \eea l.h.s. is taken as $+\iy$ if no infimum exists in the above two definitions. Note that $\sigma_{bd}$ is the first hitting time of the boundary,   while $\vartheta_{R^{-1}a}$ is the first entrance time into the open upper orthant $\{x\gg R^{-1}a\}$ for $\{w_k:1\leq k\leq n\}.$

\begin{Lemma}\label{Lem:Connc1} Let $R$ satisfy (H1). Let $n\geq 1$ be fixed, and $a\in\overline{G}.$ If $\Delta y^{(a)}_n \gg 0,$ then $\vartheta_{R^{-1}a} \leq n < \sigma_{bd},$ and $w_n \gg R^{-1}a.$
\end{Lemma}
\proof By the string \eqref{ssruinequi}, $\sum_{\ell=1}^k \hat{u}_{\ell} =\sum_{\ell=1}^k (-R^{-1}u_{n+1-\ell}) \gg 0,$ for $1\leq k\leq n,$ and $\sum_{\ell=1}^n \hat{u}_{\ell} =\sum_{\ell=1}^n (-R^{-1}u_{\ell}) \gg R^{-1}a.$ So by definition of $LCP(w_{k-1}+\hat{u}_k,R^{-1})$ and \eqref{seqhat}, we now get $\Delta v_k =0,~w_k \gg 0$ for $1\leq k\leq n,$ and $w_n \gg R^{-1}a.$ Result now follows by definitions \eqref{hittimebd},\eqref{exitss}. \qed

Our next objective is to prove a converse of Lemma~\ref{Lem:Connc1}. If $\{w_k:1\leq k\leq n\}$ does not hit $\partial G,$ and $w_n \gg R^{-1}a,$ in the phraseology of storage network, note the following. At $k=n-1,$ $w_{n-1}$ is more than sufficient to meet any potential reinforcement required due to $(\hat{u}_n-R^{-1}a),$ (in the sense that $w_{n-1}$ is enough to meet any reinforcement that may be required due to $\hat{u}_n$ and still be left with a stock of at least $R^{-1}a.$) And at $k=1,2,\cs,(n-2),$ $w_k$ is more than sufficient to meet any potential reinforcement required due to $\hat{u}_{k+1},\cs,\hat{u}_{n-1},(\hat{u}_n-R^{-1}a).$ For fixed $1\leq k\leq (n-1),$ note that part of the potential reinforcement required due to $\hat{u}_{k+m}$ can be met from $\hat{u}_{k+1},\cs,\hat{u}_{k+m-1}.$

The above comments lead us to the following finite auxiliary sequence of LCP's. Let $\Delta\xi^{(a)}_1=\Phi((-R^{-1}a-R^{-1}u_1),R^{-1}),~\z^{(a)}_1=\Psi((-R^{-1}a-R^{-1}u_1),R^{-1}),$ and $\Delta\xi^{(a)}_k=\Phi((-R^{-1}\Delta\xi^{(a)}_{k-1}-R^{-1}u_k),R^{-1}),~\z^{(a)}_k=\Psi((-R^{-1}\Delta\xi^{(a)}_{k-1}-R^{-1}u_k),R^{-1}),$ for $2\leq k\leq n.$ Therefore we have \bea\label{auxlcp1} \z^{(a)}_1 & = & -R^{-1}a-R^{-1}u_1 + R^{-1}\Delta\xi^{(a)}_1,  \\ \label{auxlcpk} \z^{(a)}_k & = & -R^{-1}\Delta\xi^{(a)}_{k-1}-R^{-1}u_k + R^{-1}\Delta\xi^{(a)}_k,~2\leq k\leq n, \eea subject to $\z^{(a)}_k\geq 0,~\Delta\xi^{(a)}_k\geq 0,~\langle\z^{(a)}_k, \Delta\xi^{(a)}_k \rangle =0,~1\leq k\leq n.$ It may be noted that the above auxiliary sequence of LCP's does not form an SP. However, we have the following.

\begin{Lemma}\label{Lem:auxlcp} Let $R$ satisfy (H1). Then $\Delta\xi^{(a)}_k = z^{(a)}_k,~\z^{(a)}_k = \Delta y^{(a)}_k$ for $1\leq k\leq n.$
\end{Lemma}
\proof For $k=1$, the result is immediate from Proposition~\ref{Prop:lcphat}. For $k\geq 2,$ by repeated use of Proposition~\ref{Prop:lcphat}, we get  \Bea \Delta\xi^{(a)}_k & = & \Phi(-R^{-1}\Delta\xi^{(a)}_{k-1}-R^{-1}u_k,R^{-1}) \\ & = & \Psi(-R(-R^{-1}\Delta\xi^{(a)}_{k-1}-R^{-1}u_k),R) \\ & = & \Psi(\Delta\xi^{(a)}_{k-1} + u_k,R) \\ & = & \Psi(z^{(a)}_{k-1}+u_k,R) = z^{(a)}_k, \Eea as required. The other assertion is similarly proved. \qed

\begin{Lemma}\label{Lem:Connc2} Let $R$ satisfy (H1). Let $n\geq 1$ be fixed, and $a\in\overline{G}.$ If $\vartheta_{R^{-1}a} \leq n < \sigma_{bd}$ and $w_n \gg R^{-1}a,$ then $\Delta y^{(a)}_n \gg 0.$
\end{Lemma}
\proof The discussion following Lemma~\ref{Lem:Connc1} indicates that we first look at $k=n-1.$ To find the potential reinforcement required due to $\hat{u}_n-R^{-1}a = (-R^{-1}u_1-R^{-1}a)$ one needs to solve $LCP((-R^{-1}a-R^{-1}u_1),R^{-1}).$ From \eqref{auxlcp1} it is clear that $R^{-1}\Delta\xi^{(a)}_1$ is the potential reinforcement required due to $(-R^{-1}a-R^{(a)}u_1).$ So by the hypothesis, it follows that $w_{n-1} \gg R^{-1}\Delta\xi^{(a)}_1.$ Proceeding analogously for $1\leq k\leq n,$ by \eqref{auxlcpk} we see that at time $k,$ the potential reinforcement required due to $\hat{u}_{k+1},\cs,\hat{u}_{n-1},(\hat{u}_n-R^{-1}a)$ is $R^{-1}\Delta\xi^{(a)}_{n-k}.$ Hence by our hypothesis it now follows that $w_k \gg R^{-1}\Delta\xi^{(a)}_{n-k},~k=n-1,n-2,\cs,2,1.$ In particular $w_1 \gg R^{-1}\Delta\xi^{(a)}_{n-1}.$ As $w_1=\hat{u}_1=-R^{-1}u_n$ by Lemma~\ref{Lem:auxlcp} it now follows that $-R^{-1}u_n \gg R^{-1}z^{(a)}_{n-1}.$ Hence it follows by Proposition~\ref{Prop:3ruin} that $\Delta y^{(a)}_n \gg 0.$      \qed

Combining Lemma~\ref{Lem:Connc1} and Lemma~\ref{Lem:Connc2} we get first part of the next result.
\begin{Theorem}\label{Th:Conncss} (i) Let $R$ satisfy (H1). Let $n\geq 1$ be fixed, and $a\in\overline{G}.$ Then $\Delta y^{(a)}_n \gg 0$ if and only if $\vartheta_{R^{-1}a} \leq n < \sigma_{bd},$ $w_n \gg R^{-1}a.$ In such a case $v_n=0,~z^{(a)}_n=0,~y^{(a)}_n=-R^{-1}a+\sum_{\ell=1}^n(-R^{-1}u_{\ell}).$

(ii) With $R,n$ as in (i), $[y^{(0)}_n : \Delta y^{(0)}_n\gg 0]   = [ w_n : n<\sigma_{bd}] $ holds; (here $[\gamma : A]$ denotes the value of $\gamma$ subject to the constraint $A.$) In such a case  \bea\label{ywconncss}  y^{(0)}_n  =  w_n  & = & \sum_{\ell=1}^n (-R^{-1}u_{\ell}). \eea
\end{Theorem}
\proof To prove part (ii), take $a=0$ in part (i). Clearly \eqref{ywconncss} holds in this case. \qed

The above discussion can be extended to quantities related to other two notions of ruin as well. For this define \bea\label{hittime0} \sigma_0 &=& \inf\{k\geq 1: w_k=0\}, \\ \label{exits} \theta_{R^{-1}a} & = & \inf\{k\geq 1: w_k > R^{-1}a \}. \eea Note that $\sigma_0$ is the first hitting time of the origin $0,$ while $\theta_{R^{-1}a}$ is the first entrance time into the set $[\{x\geq R^{-1}a\}\setminus\{R^{-1}a\}]$ for $\{w_k:1\leq k\leq n\}.$  We have

\begin{Theorem}\label{Th:Conncs} (i) Let $R$ satisfy (H1). Let $n\geq 1$ be fixed, and $a\in\overline{G}.$ Then $z^{(a)}_n=0,~\Delta y^{(a)}_n > 0$ if and only if $\theta_{R^{-1}a} \leq n < \sigma_0,~v_n=0,~$ $w_n > R^{-1}a.$ In such a case $v_n=0,~z^{(a)}_n=0,~y^{(a)}_n=-R^{-1}a+\sum_{\ell=1}^n(-R^{-1}u_{\ell}).$

(ii) With $R,n$ as in (i), $[y^{(0)}_n : z^{(0)}_n=0, \Delta y^{(0)}_n >0]   = [w_n: v_n=0, n< \sigma_0] $ holds; in such a case also \eqref{ywconncss} holds.
\end{Theorem}
\proof As $\Delta y^{(a)}_n >0$ does not necessarily imply $z^{(a)}_n=0,$ we need to specify it as well; similarly $v_n=0$ has to be spelt out. With these modifications, replacing $\gg$ by $>$ at appropriate places in the earlier discussion/ results, the theorem can be established.   \qed

For $\{w_k:1\leq k\leq n\},$ denote the first entrance time into the closed upper orthant $\{x\geq R^{-1}a\}$ by \bea\label{exit} \overline{\vartheta}_{R^{-1}a} &=& \inf\{k \geq 1:w_k \geq R^{-1}a \}. \eea Similar analysis, replacing $\gg$ by $\geq$ leads to

\begin{Theorem}\label{Th:Connc} (i) Let $R$ satisfy (H1). Let $n\geq 1$ be fixed, and $a\in\overline{G}.$ Then $z^{(a)}_n=0$ if and only if $\overline{\vartheta}_{R^{-1}a} \leq n,~v_n=0,~$ $w_n \geq R^{-1}a.$ In such a case $v_n=0,~z^{(a)}_n=0,~y^{(a)}_n=-R^{-1}a+\sum_{\ell=1}^n(-R^{-1}u_{\ell}).$

(ii) With $R,n$ as in (i), $[y^{(0)}_n : z^{(0)}_n=0]   =[ w_n : v_n=0]$ holds; in such a case also \eqref{ywconncss} holds.
\end{Theorem}

An interesting corollary of the above is
\begin{Corollary}\label{Cor:y0ya} Notation as above. Then

(i) $\Delta y^{(a)}_n \gg 0$ is equivalent to $\Delta y^{(0)}_n \gg 0,~y^{(0)}_n \gg R^{-1}a;$

(ii) $z^{(a)}_n =0$ is equivalent to $z^{(0)}_n=0,~y^{(0)}_n \geq R^{-1}a.$

In such a case $y^{(a)}_k=y^{(0)}_k-R^{-1}a,~z^{(a)}_k=z^{(0)}_k$ for $k\geq n.$
\end{Corollary}

\section{ Stochastic setup }
\setcounter{equation}{0}
In Section 2 we had derived some sample path duality results. We now consider the corresponding situation in the stochastic setup. The connection among the ruin probability of an insurance network with initial capital $a,$ certain asymptotic behaviour of the storage network starting at $0,$ and an appropriate asymptotic functional of the pushing process of the insurance network with initial capital $0,$ is our goal. We make the assumptions (H1)-(H8) at various stages.
\subsection{ Two related regulated random walks }

Assume (H1). Put $U_k(\omega)=A_k(\omega) c - X_k(\omega),~\omega\in\Omega,k=1,2,\cs.$ Let $a\in\overline{G}$ denote the vector of initial capitals. Solving the deterministic problem $SP(\{a+\sum U_{\ell}(\omega)\},R)$ path-by-path, that is, taking $u_n=U_n(\omega),~n\geq 1$ for an arbitrary but fixed $\omega\in\Omega,$ we get the pushing process $\{Y^{(a)}_n:n\geq 0\},$ and regulated process $\{Z^{(a)}_n:n\geq 0\}$ satisfying
\begin{description}
\item{\bf (S0)} $Y^{(a)}_0=0,~Z^{(a)}_0=a.$
\item{\bf (S1)} Skorokhod equation holds, that is, \bea\label{SEq} Z^{(a)}_n(\omega) & = & a + \sum_{\ell=1}^n U_{\ell}(\omega) + RY^{(a)}_n(\omega) \nonumber \\ & = & Z^{(a)}_{n-1}(\omega) + U_n(\omega) + R \Delta Y^{(a)}_n(\omega),~~n\geq 1, \eea where $\Delta Y^{(a)}_n(\omega)= Y^{(a)}_n(\omega)-Y^{(a)}_{n-1}(\omega).$
\item{\bf (S2)} $Z^{(a)}_n(\omega)\in\overline{G},~~n\geq 1.$
\item{\bf (S3)} $ \Delta Y^{(a)}_n(\omega)\geq 0$ as vectors; also \bea\label{MinCon}
\langle Z^{(a)}_n(\omega),\Delta Y^{(a)}_n(\omega)\rangle  & = & 0,~~n\geq 1. \eea
\end{description}
So the pair of processes $\{Y^{(a)}_n,Z^{(a)}_n:n\geq 0\}$ solves the Skorokhod problem $SP(\{a+\sum U_k\},R).$

Suppose (H3)-(H5) hold in addition to (H1). Then the partial sums $(a+\sum_{\ell=1}^nU_{\ell}),~n\geq 1$ form a random walk in $\R^d.$ So $Z^{(a)}_n,~n\geq 0$ is a regulated random walk starting at $a$ with $Y^{(a)}_n$ being the corresponding pushing process. This set up represents a discrete time insurance network operating under a risk diversifying treaty.

Next, put $\hat{U}_k(\omega) = - R^{-1}U_k(\omega),~\omega\in\Omega,k=1,2,\cs.$ Take $W_0=0,V_0=0.$ Solving the linear complementarity problem $LCP((W_{k-1}(\omega)+\hat{U}_k(\omega)),R^{-1}),~k=1,2,\cs$ recursively we get the pushing process $\{V_n:n\geq 0\},$ and regulated process $\{W_n:n\geq 0\}$ satisfying
\begin{description}
\item{\bf (DS0)} $V_0=0,~W_0=0.$
\item{\bf (DS1)} Skorokhod equation holds, that is, \bea\label{hatSEq} W_n(\omega) & = &  \sum_{\ell=1}^n \hat{U}_{\ell}(\omega) + R^{-1}V_n(\omega) \nonumber \\ & = & W_{n-1}(\omega) + \hat{U}_n(\omega) + R^{-1} \Delta V_n(\omega),~~n\geq 1, \eea where $\Delta V_n(\omega)= V_n(\omega)-V_{n-1}(\omega).$
\item{\bf (DS2)} $W_n(\omega)\in\overline{G},~~n\geq 1.$
\item{\bf (DS3)} $ \Delta V_n(\omega)\geq 0$ as vectors; also \bea\label{DMinCon}
\langle W_n(\omega),\Delta V_n(\omega)\rangle  & = & 0,~~n\geq 1. \eea
\end{description}
Thus the pair of processes $\{V_n,W_n:n\geq 0\}$ solves the Skorokhod problem $SP(\{\sum \hat{U}_k\},R^{-1}).$
Under the hypotheses (H1),(H3)-(H5), as before, $W_n,~n\geq 0$ is a regulated random walk in the orthant starting at $0.$ This is the stochastic analogue of the storage network considered earlier, but now over an infinite time horizon.

The next two results give implications of the coordinatewise net profit condition; see also the proof of Proposition 2.2 of \cite{Ra2012}.

\begin{Proposition}\label{Prop:Yntight} Let (H1),(H3)-(H5),(H8) hold. Then there is a $\overline{G}-$valued random variable $H_0$ such that \bea\label{Ynbd} \P(Y^{(0)}_n \ll H_0,~\forall n\geq 0) & = & 1. \eea Moreover, for a.e.$\o,$ there exists an integer $k(\o)$ such that
$(\Delta Y^{(0)}_k(\o))_i=0,1\leq i\leq d$ for all $k\geq k(\o).$ 

\end{Proposition}
\proof By (H8) note that $\E((U_1)_i)>0,~1\leq i\leq d.$ Also by (H3)-(H5), $\{(U_{\ell})_i,\ell\geq 1\}$ is a sequence of i.i.d. random variables. So by the strong law of large numbers there is $B\in\mathcal{F}$ with $\P(B)=1$ such that $\sum_{\ell=1}^n (U_{\ell}(\omega))_i \r +\iy$ as $n\r\iy,$ for all $1\leq i\leq d,$ for $\omega\in B.$ Hence for $\omega\in B$ there is an integer $n_0(\omega)$ such that
$$\sup_{k\leq n} \max \{0,-(\sum_{\ell =1}^k (U_{\ell}(\omega))_i)\} = (h_0(\omega))_i \triangleq  \sup_{k\leq n_0(\omega)} \max \{0,-(\sum_{\ell =1}^k (U_{\ell}(\omega))_i)\}, $$ for $n\geq n_0(\omega).$ Put $h_0(\omega)=((h_0(\omega))_1,\cs,(h_0(\omega))_d),$ and take $H_0(\omega)=R^{-1}h_0(\omega),~\omega\in B.$ Observe that \eqref{Ynbd} now follows from \eqref{apri}. From the above and \eqref{apri} we also get that 
$\Delta Y^{(0)}_n(\o)=0$ for $n>n_0(\o),\o\in B.$    \qed

Analogous to \eqref{hittimebd},\eqref{hittime0} respectively, define \bea\label{Hitbd}  \sigma_{bd}(\omega)  &=& \inf\{k\geq 1: W_k(\omega)\in\partial G \}, \\
\label{Hit0} \sigma_0(\omega)  &=& \inf\{k\geq 1: W_k(\omega) = 0 \}.  \eea

\begin{Proposition}\label{Prop:Hitbd} Let (H1),(H3)-(H5),(H8) hold. Then \bea\label{hitbdfin} \P(\sigma_{bd} < \iy) &=& 1. \eea
\end{Proposition}
\proof Let $B$ be as in the proof of Proposition~\ref{Prop:Yntight}.
By \eqref{Rinv}, $R^{-1}$ has only nonnegative entries with diagonal entries $\geq 1.$ Hence $\sum_{\ell=1}^n (-R^{-1}U_{\ell}(\omega))_i \r (-\iy)$ as $n\r\iy,$ for all $1\leq i\leq d,$ for $\omega\in B.$ So for $\omega\in B$ there is an integer $n_1(\omega)$ such that  $\sum_{\ell=1}^n (-R^{-1}U_{\ell}(\omega))_i < 0 $ for all $n\geq n_1(\omega),1\leq i\leq d.$ Required conclusion \eqref{hitbdfin} now follows.   \qed

\begin{Remark}\label{Rem:eqllaw} {\rm Suppose $n\geq 1$ is arbitrary but fixed. Note that the random variable $(U_n,U_{n-1},\cs,U_1)$ has the same distribution as $(U_1,\cs,U_{n-1},U_n).$ Hence the pathwise discussion in Section 2 for $w_1,w_2,\cs,w_{n-1},w_n$ is applicable to random variables $W_1,W_2,\cs,W_{n-1},W_n$ in the sense of results holding in law, that is,
equality in law; ($=^d$ shall denote equality in law.)}
\end{Remark}

\subsection{ Ruin of insurance network }

We now define ruin times corresponding to the notions of ruin introduced to earlier. For $a\in\overline{G},\omega\in\Omega$ define  \bea\label{ssruintime} \varrho^{(a)}_{ss}(\omega) & = &  \inf\{k\geq 1: \Delta Y^{(a)}_k(\omega) \gg 0\}, \\ \label{sruintime} \varrho^{(a)}_s(\omega) & = & \inf\{k\geq 1: Z^{(a)}_k(\omega) =0, \Delta Y^{(a)}_k(\omega) > 0\}, \\ \label{ruintime} \varrho^{(a)}_r(\omega) & = & \inf\{k\geq 1: Z^{(a)}_k(\omega) =0 \}. \eea 
$\varrho^{(a)}_{ss}(\o)$ is taken to be $+\iy$ if there is no $k\geq 1$ satisfying the requirement; similar comment applies to the other cases. Clearly $\varrho^{(a)}_r \leq \varrho^{(a)}_s \leq \varrho^{(a)}_{ss}.$ Note that $\varrho^{(a)}_{ss}, \varrho^{(a)}_{s}, \varrho^{(a)}_r$ are ruin times corresponding to, respectively, ss-ruin, s-ruin, ruin.

The next result indicates when the three notions may coincide.

\begin{Proposition}\label{Prop:equiruin} Let (H1),(H3)-(H5),(H7) hold, and $c \gg 0.$ Then for any $a\in\overline{G},$ \bea\label{equiruin} \P(\varrho^{(a)}_r = \varrho^{(a)}_s = \varrho^{(a)}_{ss}) & = & 1. \eea
\end{Proposition}
\proof For fixed integer $k\geq 1,$ we need to show that $Z^{(a)}_k=0$ implies $\Delta Y^{(a)}_k \gg 0$ with probability $1.$ As $c \gg 0$ and $A_{\ell} >0,$ note that $(R^{-1}(Z^{(a)}_{k-1}+A_k c))_i>0$ for all $1\leq i\leq d.$ By (H7) $(R^{-1}X_k)_i\geq 0,$ and has no atoms on $(0,\iy).$ Also by (H3)-(H5), $R^{-1}(Z^{(a)}_{k-1}+A_k c)$ and $R^{-1}X_k$ are independent random variables. Consequently $\P[(R^{-1}(Z^{(a)}_{k-1}+A_k c -X_k))_i=0]=0$ for all $1\leq i\leq d.$ If $Z^{(a)}_k(\omega)=0,$ by \eqref{SEq} note that $-(R^{-1}(Z^{(a)}_{k-1}+A_k c -X_k))_i(\omega) = (\Delta Y^{(a)}_k)_i(\omega) \geq 0.$ The required conclusion now follows. \qed

The next result implies that various events associated with ss-ruin, like stochastic analogues of the string  \eqref{ssruinequi}, have positive probability.

\begin{Proposition}\label{Prop:nontrivruin} Let (H1)-(H6) hold, and let $c\gg 0.$ Then the following hold.

(i) $\delta_0=\P(-R^{-1}U_{\ell} \gg 0) > 0  $ for any $\ell\geq 1.$ Moreover $\P( -R^{-1}U_{\ell} \gg 0~~{\rm infinitely}~{\rm often}~)=1.$

(ii) $\P(-R^{-1}U_{\ell} \gg b) >0$ for any fixed $\ell\geq 1,~b\in\overline{G}.$
\end{Proposition}
\proof Let $\ell\geq 1.$ Let $k$ be as in (H2). Take $\beta = \max\{(R^{-1}c)_j/(R^{-1})_{jk}:1\leq j\leq d\}.$ By (H2) we get $0<\beta<\iy.$ By (H6) support of $(X_{\ell})_k = [0,\iy).$ Consequently, as $A_{\ell}>0$ and independent of $(X_{\ell})_k,$ we now get $\P((X_{\ell})_k > \beta A_\ell ) >0.$ So by the definition of $\beta$ note that $\P(-(R^{-1}U_{\ell})_i >0, \forall 1\leq i\leq d ) > 0.$ The first assertion now follows. As $U_{\ell}$ are i.i.d. random variables, note that $\delta_0$ does not depend on $\ell.$ An application of the second Borel-Cantelli lemma now gives the second assertion. This proves (i).

To prove (ii), let $\beta_0=\max\{(b)_j/(R^{-1})_{jk}:1\leq j\leq d\}.$  Proceeding as in (i), we get that  $\P((X_{\ell})_k > \beta A_\ell +\beta_0) >0.$ Required conclusion is obtained as above. \qed

\subsection{ Pollaczek-Khinchine formula } 

We will now confine ourselves to the ss-ruin problem. Besides being the appropriate 
$d-$dimensional analogue of the classical ruin problem, it seems to be more amenable 
to analysis. Of course, Proposition~\ref{Prop:equiruin} gives sufficient conditions 
for the three notions of ruin to be equivalent with probability one. We begin with a 
duality result. 

\begin{Theorem}\label{Th:duality1} Assume (H1)-(H6),(H8). Let $a\in\overline{G}.$ Then
\bea\label{ssduality} 0 < \P( \varrho^{(a)}_{ss} < \iy) &=&
\P(\vartheta_{R^{-1}a} < \sigma_{bd}) <1. \eea where
$\vartheta_{R^{-1}a}(\omega) = \inf\{k\geq 1: W_k(\omega) \gg
R^{-1}a\},$ and $\varrho^{(a)}_{ss},\sigma_{bd}$ are given
respectively by \eqref{ssruintime},\eqref{Hitbd}. Moreover
$\P(\Delta Y^{(a)}_n \gg 0)>0$ and hence $\P(\sigma_{bd} >n)>0$ for
any $n\geq 1.$ Also \bea\label{survivprob} \lim_{|a|\r\iy,a\in G}
\P( \varrho^{(a)}_{ss} < \iy) &=& 0. \eea
\end{Theorem}
\proof Note that the equality in \eqref{ssduality} follows by part (i) of Theorem~\ref{Th:Conncss}. Now taking $b=R^{-1}a$ in part (ii) of Proposition~\ref{Prop:nontrivruin}, we see that $\P(\varrho^{(a)}_{ss}=1)>0;$ so the first inequality in \eqref{ssduality} follows. By (H8), $\E((R^{-1}U_1)_i)>0,$ and hence $\P((R^{-1}U_1)_i>0)>0$ for all $i.$ Hence $\P(W_1\in\partial G)>0;$ that is, $\P(\sigma_{bd}=1)>0.$ Consequently the second inequality in \eqref{ssduality} now follows. Next fix $n\geq 1;$ we proceed as in the proof of part (i) of Proposition~\ref{Prop:nontrivruin}.  Note that $R^{-1}X_n$ is independent of $R^{-1}Z^{(a)}+A_nR^{-1}c;$ so by (H2),(H6) we get $\P(R^{-1}X_n \gg R^{-1}Z^{(a)}_{n-1}+A_nR^{-1}c) >0;$ that is, $\P(-R^{-1}U_n \gg R^{-1}Z^{(a)}_{n-1}) >0,$ which is equivalent to ss-ruin Proposition~\ref{Prop:3ruin}. So by part (i) of Theorem~\ref{Th:Conncss}, the second assertion follows. Finally, by Proposition~\ref{Prop:Hitbd}, it follows that $\lim_{n\r\iy} \P(\sigma_{bd}>n)=0.$ Hence \eqref{ssduality} now implies \eqref{survivprob}.   \qed

Now assume (H1)-(H6),(H8). Note that the $d-$dimensional
insurance network $\{(Y^{(0)}_n,Z^{(0)}_n):n=0,1,2,\cs\},$ with initial capital $0,$
is a strong Markov process
starting at $(0,0).$ As there is no dispersion, and as drift, reflection are constants,
$\{Z^{(0)}_n:n=0,1,2,\cs\}$ is also a strong Markov process starting at $0.$ Take
$\t_0\equiv 0.$ For $n\geq 1,$ define
\bea\label{taun} \t_n(\o) &=& \inf\{k\geq \t_{n-1}(\o)+1:\Delta Y^{(0)}_k(\o)\gg 0\}, \eea
if the set $\{k\geq \t_{n-1}(\o)+1:\Delta Y^{(0)}_k(\o)\gg 0\} \neq \varnothing;$
put $\t_n(\o)=+\iy$ if there is no $k\geq\t_{n-1}(\o)+1$ such that
$\Delta Y^{(0)}_k(\o)\gg 0.$  Note that $\t_k,k\geq 0$ are stopping times w.r.t.
the natural filtration.

For convenience, write $Y^{(0)}(k,\o)=Y^{(0)}_k(\o),Z^{(0)}(k,\o)=Z^{(0)}_k(\o);$
so $Y^{(0)}(\t_j,\o)=Y^{(0)}(\t_j(\o),\o),Z^{(0)}(\t_j,\o)=Z^{(0)}(\t_j(\o),\o).$ For $n\geq 1,$
define \bea\label{ladhtn} L_n(\o) &=& Y^{(0)}(\t_n,\o)-Y^{(0)}(\t_{n-1},\o),~~{\rm if}~~\t_n(\o)<\iy, \nonumber \\
&=& 0,~~{\rm if}~~\t_n(\o)=+\iy; \\
L^+_n(\c) &=& L_n(\c)~~{\rm restricted}~{\rm to}~~\{\t_n<\iy\}; \eea in the above note that
$Y^{(0)}(\t_0)\equiv 0.$ Clearly $L_n$ takes value in $\{0\}\cup G,$ and $L^+_n$ in $G.$ We shall call $L^+_1$
the \emph{$d-$dimensional first strictly ascending ladder height} random variable, and $L^+_k$
the \emph{$d-$dimensional $k-$th strictly ascending ladder height} random variable.

Now by the second assertion in Proposition~\ref{Prop:Yntight}, note that  
for a.e.$\o\in\Omega,$ there is $n_0(\o)$ such that $\t_n(\o)=+\iy,$
and hence $L_n(\o)=0$ for all $n\geq n_0(\o).$ Define
\bea\label{beta} \beta (\o) &=& \inf\{k\geq 1:\t_k(\o)=+\iy\}=
\inf\{k\geq 1:L_k(\o)=0\}. \eea From the above it is clear that $\b
<\iy$ with probability $1.$

Note that $Z^{(0)}(\t_n,\o)=0$ if $\t_n(\o)<\iy.$ Consequently, by the strong Markov
property, conditional on $\t_n<\iy,$ $$\{((Y^{(0)}(\t_k+j)-Y^{(0)}(t_k)),Z^{(0)}(\t_k+j)):
j=0,1,2,\cs \},~~k=0,1,2,\cs,n$$ are $(n+1)$ independent stochastic processes; also the
first $n$ of these, that is, corresponding to $k=0,1,\cs,n-1$ are identically distributed
as well. In particular, conditional on $\t_n<\iy,$ while the $(n+1)$
random variables $\t_1,\t_2-\t_1,\cs,\t_n-\t_{n-1},\t_{n+1}-\t_n$ are independent, the first
$n$ of these are i.i.d. finite valued random variables; (here $\t_{n+1}-\t_n$ is taken as
$+\iy$ when $\t_{n+1}=+\iy.$) Hence, conditional on $\t_n<\iy,$ random vectors $L_1,L_2,
\cs,L_n,L_{n+1}$ are independent, and $L_1,\cs,L_n$ are i.i.d. random vectors taking value
in $G,$ with $L_k=L^+_k,1\leq k\leq n;$ also $L_{n+1}=0$ if and only if $\t_{n+1}=+\iy.$ For any  Borel set $B\subseteq G$ define
\bea\label{ladhtddim} \a_+(B) &=&\P(L^+_1\in B),~~B\subseteq G, \\
\a_0(B) &=& \frac{1}{\a_+(G)}\a_+(B),~~B\subseteq G. \eea Note that
$\a_+$ is a defective distribution, while $\a_0$ is the
corresponding normalized probability distribution both concentrated
on $G.$ Take $M_0\equiv 0,$ and for $n\geq 1,$ define
\bea\label{MnM} M_n(\o) &=& \sum_{j=1}^nL_j(\o), \nonumber \\
M(\o) & = & \sum_{j=1}^{\iy}L_j(\o); \eea note that
$M_n(\o)=Y^{(0)}(\t_n(\o),\o),$ if $\t_n(\o)<\iy$ and $M_n(\o)=M_{n-1}(\o)$ if $\t_n(\o)=\iy.$

\begin{Theorem}\label{Th:PKM} Assume (H1)-(H6),(H8). Let $\varrho^{(a)}_{ss},\sigma_{bd},\vartheta_{R^{-1}a}$
be as in Theorem~\ref{Th:duality1}. Denote $p \triangleq
\P(\hat{U}_1\in G)=\P(-R^{-1}U_1\in G);$ note that $0<p<1.$ Then
$(\b-1)$ has a geometric distribution with parameter $(1-p),$
$\a_+(G)=p,$ and the distribution of $M$ is the geometric compound
\bea\label{PKM} \nu_M(B) & = & (1-p)\delta_0(B) + \sum_{k=1}^{\iy}
(1-p)p^k\a_0^{\ast(k)}(B),  \eea  for any Borel set $B\subseteq
\{0\} \cup G.$ Moreover ruin probability for the insurance network
is given by
\bea\label{PKruin} \P( \varrho^{(a)}_{ss} < \iy) &=& \P(M \gg R^{-1}a ) \nonumber \\
&=& (1-p) \sum_{n=1}^{\iy} \a_+^{\ast(n)} (\{x \gg R^{-1}a \}),~~a\in\overline{G}. \eea
\end{Theorem}

\proof By hypotheses (H5),(H6),(H1) and \eqref{Rinv} we get $p>0;$
by the coordinatewise net profit condition (H8) it is clear that
$p<1.$ As $0=R^{-1}0,$ we have $\vartheta_0=\inf\{k\geq 1:W_k\in
G\}.$ Clearly $\{\sigma_{bd} \leq \vartheta_0\} =\{\sigma_{bd} <
\vartheta_0\} =\{\sigma_{bd}=1\}$ as events. Consequently by
Theorem~\ref{Th:duality1} \bea\label{beta1} \P(\b =1) &=&
\P(\t_1=+\iy) = \P(\Delta Y^{(0)}_k \gg 0~{\rm never}~ {\rm
happens}) \nonumber \\ &=& \P(\varrho^{(0)}_{ss}=\iy) =
\P(\sigma_{bd} < \vartheta_0) \nonumber \\ &=& \P(\sigma_{bd} = 1) =
(1-p). \eea From the discussion preceding the theorem, conditional
on $\t_n<\iy,$ note that $L_1,\cs,L_n$ are i.i.d. random vectors
with distribution $\a_0$
for each $n;$ also \bea\label{L=0} \P(L_{n+1}=0\mid \t_n<\iy)&=&\P(\t_{n+1}=\iy\mid \t_n<\iy) \nonumber \\
& = & \P(\t_1 = \iy) =\P(L_1 =0) = (1-p). \eea Since
$\{\t_j<\iy,j\leq k\}=\{\t_k<\iy\},k\geq 1,$ proceeding recursively
and using \eqref{beta1},\eqref{L=0}, we have
\bea\label{betan} \P(\b = n+1) &=& \P(\t_{n+1}=\iy,\t_j<\iy,1\leq j\leq n) \nonumber \\
&=& \P(\t_{n+1}=\iy\mid \t_n<\iy)\P(\t_n<\iy) = (1-p)\P(\t_n<\iy) \nonumber \\
&=& (1-p)\P(\t_n<\iy\mid \t_{n-1}<\iy)\P(\t_{n-1}<\iy) \nonumber \\
&=& (1-p)\P(\t_1<\iy)\P(\t_{n-1}<\iy)\nonumber \\ & = &(1-p)p\P(\t_{n-1}<\iy)
= (1-p)p^n,~~n=0,1,\cs \eea Thus $\b-1$ has a geometric distribution with parameter $1-p.$
From \eqref{L=0} it also follows that $\a_+(G)=p.$

Clearly \bea\label{Mbeta} M(\o) &=& \lim_{n\r\iy} M_n(\o) =\sup_{n\geq 0} M_n(\o) \nonumber \\
&=&  \sum_{k=1}^{\b(\o)-1} L_k(\o),~~{\rm if}~~\b(\o)\geq 2,  \\
&=& 0,~~{\rm if}~~\b(\o)=1. \nonumber \eea By \eqref{Mbeta},\eqref{beta1}, we get
$\P(M=0)=(1-p).$ Now by \eqref{MnM}, conditional on $\b =n+1,$ it is seen that $M_n$ is
distributed as $\a_0^{\ast(n)}.$ Hence, using \eqref{betan},\eqref{ladhtddim}, for any Borel set $B\subseteq G,$
\bea\label{Mdist} \P(M\in B) &=& \sum_{n=1}^{\iy} \P(M\in B\mid \b -1=n)\P(\b -1=n) \nonumber \\
& = & \sum_{n=1}^{\iy} \P(M_n\in B\mid \b=n+1)\P(\b=n+1) \nonumber \\
& = & \sum_{n=1}^{\iy} (1-p)p^n \a_0^{\ast(n)}(B) \nonumber \\
& = & (1-p) \sum_{n=1}^{\iy} \a_+^{\ast(n)}(B). \eea

By \eqref{MnM}, note that
\bea\label{my0} M(\o) &=& \lim \{Y^{(0)}_k(\o): (\Delta Y^{(0)}_k(\o) \gg 0 \}. \eea Therefore
by the definition of ruin, Corollary~\ref{Cor:y0ya}, Theorem~\ref{Th:duality1} we get for any $a\in \overline{G},$
\bea\label{ruinaM} \P( \varrho^{(a)}_{ss} < \iy) &=& \P(\Delta Y^{(a)}_k \gg 0~~{\rm for}~{\rm some}~k\geq 1 ) \nonumber \\
&=& \P (\Delta Y^{(0)}_k \gg 0,~Y^{(0)}_k \gg R^{-1}a ~~{\rm for}~{\rm some}~k\geq 1 ) \nonumber \\ &=& \P(M \gg R^{-1}a ) \eea
From \eqref{Mdist},\eqref{ruinaM}, required conclusions \eqref{PKM},\eqref{PKruin} follow.  \qed

\begin{Theorem}\label{Th:duality3} Assume (H1)-(H6),(H8). Let $\varrho^{(a)}_{ss},\sigma_{bd},\vartheta_{R^{-1}a},p,\b,\a_+,\a_0,M,\nu_M$
be as in Theorem~\ref{Th:PKM}. Define \bea\label{hatUplus}
\hat{U}_k^+ &=& \hat{U}_k~~{\rm restricted}~{\rm to}~~\{\hat{U}_k\in
G\}, \\
\label{muplus} \mu_+(B) & = & \P(\hat{U}_1^+ \in B)= \P(\hat{U}_1
\in B),~~B\subseteq G. \eea Note that $\mu_+$ is a defective
distribution concentrated on $G,$ with $0< p = \mu_+(G) <1.$ Set
$\mu_0(\c)=\frac{1}{p}\mu_+(\c).$ Define the compound geometric
$\nu(\c)$ by \bea\label{PK1} \nu(B) & = & (1-p)\delta_0(B) +
\sum_{k=1}^{\iy} (1-p)p^k\mu_0^{\ast(k)}(B),~B\in\mathcal{B}(\R^d).
\eea Then the following hold:

(i) $(\sigma_{bd}-1)$ has a geometric distribution with parameter
$(1-p),$ and hence $\sigma_{bd}=^d\b.$

(ii) $\nu$ is a probability measure concentrated on $\{0\}\cup G,$
such that \bea\label{distnWsigmabd} \P((\max_{k<\sigma_{bd}}W_k) \in
B) & = & \P( W(\sigma_{bd}-1) \in B) \nonumber \\ & = & \nu(B),~~
B\subseteq (\{0\}\cup G); \eea here $(\max_{k <
\sigma_{bd}}W_k)(\omega) = (\max_{k <
\sigma_{bd}(\omega)}(W_k(\omega))_1,\cs,\max_{ k <
\sigma_{bd}(\omega)}(W_k(\omega))_d).$ Also \bea\label{partiallimWn}
\nu(\{x \gg  z\}) &=& \lim_{n\r\iy}\P(W_n \gg z,\sigma_{bd}
> n),~z\in\overline{G}; \eea that is, on $[0,\sigma_{bd}),$ $W_n$ 
converges in distribution to $W(\sigma_{bd}-1).$  

(iii) $M=^d \max_{k<\sigma_{bd}}W_k,$ and hence $\nu_M=\nu.$

(iv) $\mu_0=\a_0,~\hat{U}^+_1=^d L^+_1;$ that is, $\mu_+$ is the
\emph{$d-$dimensional ladder height distribution,} in other words
\bea\label{ladhtd} \P(L^+_1\in B) &=& \P(-R^{-1}(cA_1-X_1)\in
B),~B\subseteq G. \eea

(v) For any $a\in \overline{G},$ \bea\label{PK2}
\P(\varrho^{(a)}_{ss} <\iy) &=& \P((\max_{k<\sigma_{bd}}W_k) \gg
R^{-1}a)
\nonumber \\
 &=& \nu (\{x\gg R^{-1}a\}) \nonumber \\
& = & \sum_{k=1}^{\iy} (1-p) \mu_+^{\ast(k)}(\{x\gg R^{-1}a\}). \eea
\end{Theorem}

\proof We already have $0<p<1;$ also $\mu_+(G^c)=0,$ by \eqref{muplus}. Thus $\mu_+$ is a
defective distribution and $\mu_0$ a probability distribution, both
concentrated on $G.$ For $k=1,2,\cs,$ clearly $\mu_0^{\ast(k)}$ is
concentrated on $G;$ consequently the geometric compound $\nu$ is
a probability measure concentrated on $\{0\}\cup G.$

By \eqref{Hitbd} and Proposition~\ref{Prop:Hitbd}, $1\leq\sigma_{bd}
<\iy$ with probability $1.$ Also by \eqref{beta1},
$\P(\sigma_{bd}-1=0)=\P(\b-1=0)=1-p.$

By Theorem~\ref{Th:Conncss} and Remark~\ref{Rem:eqllaw}, we get for
$k=1,2,\cs $ \bea\label{eqllaw}
[W_k:\{ k<\sigma_{bd}\}]&=^d& [Y^{(0)}_k:\{\Delta Y^{(0)}_k \gg 0\}], \eea
where $[\gamma:A]$ stands for the random variable $\gamma$ restricted to the set $A.$
As $[Y^{(0)}_j:\{\Delta Y^{(0)}_j \gg 0\}] \gg Y^{(0)}_{j-1}~\forall j,$ 
$\{\sigma_{bd}>k\}=\cap_{\ell=1}^k\{\sigma_{bd}>\ell\},$ and 
$W_j=\sum_{\ell=1}^j\hat{U}_{\ell},~1\leq j\leq k$ on $\{\sigma_{bd}>k\},$
 by \eqref{eqllaw}
\bea\label{sigmaDY} \P(\sigma_{bd}>k) &=& \P(W_k \gg W_{k-1}\gg \cs \gg W_1 \gg 0)
\nonumber \\
&=& \P(\hat{U}_j \gg 0,~1\leq j\leq k)
= p^k. \eea Thus $\sigma_{bd}-1$ has a geometric distribution
with parameter $(1-p).$

Clearly $W(\sigma_{bd}-1)=W(0)=0$ if $\sigma_{bd}-1=0.$
By the arguments given in the derivation of \eqref{sigmaDY},
and \eqref{hatUplus} it follows that on
$\{\sigma_{bd}-1=k\},$ \bea\label{WsigmaxW} W(\sigma_{bd}-1) &=^d&
W_k \nonumber \\
&=^d& \sum_{j=1}^k\hat{U}_j^+ \nonumber \\
&=^d& (\max_{j< \sigma_{bd}}W_j), \eea for $k=1,2,\cs$

From \eqref{sigmaDY},\eqref{WsigmaxW}, for any Borel set
$B\subseteq (\{0\} \cup G),$ we get
\bea\label{Wsigmabdnu} \P((\max_{j<\sigma_{bd}}W_j)\in B)
&=& \P(W(\sigma_{bd}-1)\in B) \nonumber \\
&=& \sum_{k=0}^{\iy}\P(W(\sigma_{bd}-1)\in B \mid \sigma_{bd}-1=k)
\P(\sigma_{bd}-1=k) \nonumber \\
&=& (1-p)\delta_0(B) + \sum_{k=1}^{\iy}(1-p)p^k
\P(W_k\in B \mid \sigma_{bd}-1=k) \nonumber \\
&=& (1-p)\delta_0(B) + \sum_{k=1}^{\iy} (1-p)p^k \mu_0^{\ast (k)} (B) \nonumber \\
& = & \nu (B) \eea Note that \eqref{partiallimWn} is clear from the
above arguments.

To prove $M=^d \max_{k<\sigma_{bd}}W_k,$ it is enough to
consider the case $\sigma_{bd}>1,$ equivalently $\b>1.$
Now by Theorem~\ref{Th:Conncss} and
Remark~\ref{Rem:eqllaw}, we have
$$[W_k:\{ k<\sigma_{bd}\}]=^d [Y^{(0)}_k:\{\Delta Y^{(0)}_k \gg 0\}]\leq M_k,$$
where $[\gamma:A]$ stands for the random variable $\gamma$ restricted to the set $A.$
Hence by \eqref{my0}, $$\P((\max_{k<\sigma_{bd}}W_k) \gg x) \leq \P(M\gg x),~x\in\overline{G}.$$
Next, since $\b<\iy$ with probability $1,$ for a.e.$\o$ we have
$M(\o)=Y^{(0)}_k(\o)$ with $\Delta Y^{(0)}_k(\o) \gg 0,$ for
some $k$ that may depend on $\o.$ So again by Theorem~\ref{Th:Conncss} and
Remark~\ref{Rem:eqllaw} as above, $M=^dW_k$ with $k<\sigma_{bd},$
when $M=Y^{(0)}_k.$ Hence it follows that
$$\P(M\gg x)\leq \P((\max_{k<\sigma_{bd}}W_k) \gg x),~x\in\overline{G}.$$
Thus $M=^d \max_{k<\sigma_{bd}}W_k.$

Therefore Theorem~\ref{Th:PKM} and \eqref{Wsigmabdnu} imply
$\nu=\nu_M.$ As $\b=^d\sigma_{bd},$ in view of the expressions
\eqref{PKM},\eqref{PK1}, an elementary argument using characteristic
functions gives $\a_0=\mu_0,\a_+=\mu_+, L^+_1=^d\hat{U}^+_1.$

Finally \eqref{PK2} in assertion (v) follows by \eqref{PKruin},
\eqref{ruinaM}, \eqref{PK1} and assertion (iv). This completes the
proof.    \qed

\begin{Remark}\label{Rem:ladht1} {\rm In the classical one-dimensional
\emph{renewal risk (Sparre Anderson) model,} recall that "ruin" is
defined as the event that the surplus goes strictly below zero level
in finite time. This is the same as $\Delta y_n \gg 0$ for some
$n\geq 1$ in our frame work in terms of the one-dimensional
Skorokhod problem with $G=(0,\iy)$ and normal reflection; so
$R=R^{-1}=1.$ Also this is true even at the sample path level, that
is, in the deterministic set up. In the classical model,
Pollaczek-Khinchine formula for ruin problem is generally expressed
as a compound geometric distribution involving the one-dimensional
ladder height distribution; for definition of ladder height
distribution in the classical model (without any reference to
Skorokhod problem), relevant proofs and more information, see
Chapter 6 of \cite{RSST1999}. In view of Theorems~\ref{Th:PKM},
\ref{Th:duality3}, it is clear that our definition and the classical
notion of ladder height distribution coincide in the one-dimensional
case. So $\a_+$ given by \eqref{ladhtddim}, or equivalently $\mu_+$
given by \eqref{muplus}, can be regarded as the $d-$dimensional
analogue of ladder height distribution; moreover $\mu_+$ gives an
explicit expression for the ladder height distribution in the
$d-$dimensional renewal risk set up. In fact, when $d=1,$ and
$A_1,X_1$ both have exponential distributions, (that is, in the
classical Cramer-Lundberg model with exponential claim sizes,) it is
easy to verify that r.h.s. of \eqref{muplus} (or rh.s. of
\eqref{ladhtd}) is an appropriate multiple of the integrated tail
distribution of claim sizes. For the general one dimensional
Cramer-Lundberg model, it is known that the ladder height
distribution is the same as the integrated tail of claim sizes; see
\cite{RSST1999} }
\end{Remark}

{\bf Note:} \eqref{PKM} (or equivalently \eqref{PKruin}),
\eqref{PK1} (or equivalently \eqref{PK2}) may be considered
\emph{Pollaczek-Khinchine formula} for multidimensional ruin
problem. Because of \eqref{muplus}, all the quantities on r.h.s. of
\eqref{PK1}, \eqref{PK2} are in terms of the given data of the
model.

Following corollary is a version of the duality theorem in one
dimension; see \cite{AA2010}.

\begin{Corollary}\label{Cor:limdist} Let $d=1;$ define 
$\sigma_0(\o) =\inf\{k\geq 1:W_k(\o)=0\},$
and $Y^{(0)}_{\iy}(\omega) = \lim_{n\r\iy}Y^{(0)}_n(\omega).$ Then $Y^{(0)}_{\iy} = M,$
$W_n$ converges in distribution to $W(\sigma_0-1),$ and 
hence to $Y^{(0)}_{\iy};$ also for any $a\geq 0,$
\bea\label{limdist}  \P( \varrho^{(a)}_{ss} < \iy) &=& \P(W(\sigma_0-1) \gg a ) \nonumber \\
& =& \lim_{n\r\iy} \P(W_n \gg a) \nonumber \\ & = & \P(Y^{(0)}_{\iy} \gg a ). \eea
\end{Corollary}

\proof As $d=1$ note that $\partial G=\{0\},$ and hence by
Proposition~\ref{Prop:Hitbd}, $\sigma_{bd}=\sigma_0 <\iy$ with
probability $1.$ Also $Y^{(0)}_{\iy}=M,$ as $\Delta Y^{(0)}_j=0$
whenever it is not strictly positive. Clearly $W(\sigma_0)=0.$ Now
put $\gamma_0=0,\gamma_1=\sigma_0,$ and $\gamma_n=\inf\{k\geq
\gamma_{n-1}+1:W_k=0\},~n\geq 1,$ denoting successive times of visit
to the origin. As $W(\gamma_j)=0$ for all $j\geq 1,$ by the strong
Markov property of $\{W_{\ell}\},$ it follows that 
between successive visits to the origin, the process behaves like
independent copies of $\{W_0,W_1,\cs W(\sigma_0-1)\}.$ Therefore by the proof of
Theorem~\ref{Th:duality3}, we now get $W_n$ converges in
distribution to $W(\sigma_0-1);$ \eqref{limdist} is also a
consequence of Theorem~\ref{Th:duality3}. \qed

One-dimensional duality theorem is given in Corollary 3.2 on p.49 of
\cite{AA2010}, or Theorem 5.1.2 on p.151 of \cite{RSST1999}. Note
that $W_k=^d \max_{j\leq k}W_j$ if $k<\sigma_0,$ one-dimensional
analogue of a portion of \eqref{WsigmaxW}, is implicit in the proofs
given in these references, incidentally, using the one-dimensional
Skorokhod reflection map.

{\bf Note:} For $d\geq 2$ note that $M =0$ does not imply
$Y^{(0)}_n(\omega)=0$ for all $n.$ However, for $d=1$ the implication holds.

{\bf Note:} Let $d\geq 2.$ Note that $\P(\sigma_{bd} = \sigma_0) > 0;$ because
of the coordinatewise net profit condition, this probability could even be substantial.
It is not clear what conditions ensure $\P(\sigma_{bd} = \sigma_0) = 1.$ It is also not
clear when the process $\{W_n\}$ has a limiting distribution.

\begin{Remark}\label{Rem:reH2} {\rm Observe that the hypotheses (H2),
(H6) ensure that the various events associated with ruin have
positive probability even if claim size vector has only one nonzero
component, provided that the nonzero component is sufficiently
large, as proved in Proposition~\ref{Prop:nontrivruin}. In fact,
these two conditions can be replaced with
\begin{description}
\item{\textbf {(H9)}} For each $\ell\geq 1,$ $\P(R^{-1}X_{\ell} \gg x) >0$
for any $x\gg 0;$ that is, $R^{-1}X_{\ell}$ is supported on an unbounded upper orthant.
\end{description}
If claim size vector itself has an unbounded upper orthant as
support then (H9) holds as $(R^{-1})_{ii}\geq 1,\forall i$ by
\eqref{Rinv}. Note that
Theorems~\ref{Th:duality1},~\ref{Th:PKM},~\ref{Th:duality3} continue
to hold even with (H2),(H6) replaced by (H9). However, in our
exposition we persist with (H2),(H6) as we want to emphasize that
Example~\ref{Ex:renew} below is covered by our analysis. }
\end{Remark}

\begin{Remark}\label{Rem:atom0} {\rm \cite{K1996,K1997} consider storage networks
driven by fairly general Levy inputs. When the basic driving process
is of same kind, it is interesting to note the similarity and the
difference between the storage networks considered in these papers
and the one considered here. Note that $\E(\hat{U}_1) \ll  0$ by
(H8); its analogue is assumed in \cite{K1997}. While the reflection
matrix is $R$ in \cite{K1997}, in our case it is $R^{-1}$ for the
storage network. Thus in the storage network considered in
\cite{K1997} both the drift and the reflection vectors are inward
looking, while in our storage network the drift vector is inward
looking, but the reflection vector on each face of $\partial G$ is
outward looking. For  the storage network considered in \cite{K1997}
the limiting distribution has an  atom at the origin; while it is
not clear if our storage network $W_n$ has a limiting distribution
in general, note that $W(\sigma_{bd}-1)$ has an atom at $0$ in our
set up.}
\end{Remark}

\subsection{ Examples }
We give a few classes of examples for which our analysis can be applied.
\begin{Example}\label{Ex:renew} Renewal risk type network: {\rm All processes
are defined on a filtered probability space
$(\Omega,\mathcal{F},\{\mathcal{F}_t,t\geq 0\},\P).$ Let
$A_{\ell},\ell\geq 1,$ denote i.i.d. interarrival times
corresponding to a renewal counting process $\{N(t):t\geq 0\}.$ Let
$p_i,~1\leq i\leq d$ be positive numbers such that
$\sum_{i=1}^dp_i=1.$ Let $J^{(1)},\cs,J^{(d)}$ be independent one
dimensional random variables taking values in $[0,\iy);$ these are
independent of $\{A_{\ell}\}.$ Here $J^{(i)}$ represents claim size
distribution for Company $i.$ Let $J$ be a $d-$dimensional random
variable such that $J=(0,\cs,0,J^{(i)},0,\cs,0)$ with probability
$p_i,$ for $1\leq i\leq d;$ (here coordinates other than $i-$th are
zero). So $J$ takes values on the boundary of the orthant $G.$
Clearly $J$ is not absolutely continuous even if $J^{(i)}$ are.
Moreover, the marginals $(J)_i,1\leq i\leq d$ of $J$ are not
independent even though $J^{(i)}$ are. Let $F^{(i)}$ denote the
distribution function of $J^{(i)},$ and let $(F)_i$ denote the
$i-$th marginal distribution function of $J,$ for $1\leq i\leq d.$
It can be seen that $$(F)_i(u)=[p_iF^{(i)}(u) +
(1-p_i)]I_{[0,\iy)}(u),~u\in\R.$$ So, $(F)_i\neq  F^{(i)},$ and
$(F)_i$ has an atom at $u=0$ even if $F^{(i)}$ is continuous.

Now let $X_{\ell},\ell\geq 1,$ be i.i.d. random variables having
the same distribution as $J;$ these represent vector claim sizes. If
$J^{(i)}$ are continuous, note that (H7) is satisfied. Let
\bea\label{Hcts} H^{(a)}(t) & = & a + tc
+\sum_{\ell=1}^{N(t)}X_{\ell},~t\geq 0, \eea denote the joint
dynamics of the $d$ insurance companies in the absence of the risk
diversifying treaty; here $a$ denotes the initial capital of the
companies, while $c$ denotes the vector of premium rates. At an
arrival time, it is assumed that an independent mechanism governed
by the probability vector $(p_1,\cs,p_d)$ determines the company $i$
that has to take the claim; interarrival times $\{A_{\ell}\},$ claim
sizes having law $J^{(i)},1\leq i\leq d,$ and the random mechanism
governed by $(p_1,\cs,p_d)$ are taken to be independent of each
other. Thus $H^{(a)}(\c)$ is a $d-$dimensional renewal risk process.

An important special case is that of Cramer-Lundberg type network.
In the absence of treaty, the joint dynamics is that of $d$
independent Cramer-Lundberg processes $H^{(i)},1\leq i\leq d,$ with
respective initial capital $a_i,$ premium rate $(c)_i,$ claim number
process $\{N^{(i)}(t)\},$ which is a Poisson process with rate
$\lambda_i>0,$ and claim size $J^{(i)}.$ In such a case
$N(t)=N^{(1)}(t)+\cs+N^{(d)}(t),t\geq 0$ is also a Poisson process
with rate $\lambda=\lambda_1+\cs+\lambda_d;$ so $A_{\ell}$ are
i.i.d. random variables having exponential distribution with
parameter $\lambda.$ Take $p_i=\lambda_i/\lambda $ in the above set
up. Then it can be seen that $H^{(a)}(t)=^d
(H^{(1)}(t),\cs,H^{(d)}(t))$ as processes.

To describe the joint dynamics under risk diversifying treaty, 
one can use continuous time Skorokhod problem
$SP(H^{(a)}(\c),R);$ see \cite{Ra2006,Ra2012}. Accordingly, we seek
$d-$dimensional r.c.l.l. processes $\{Y^{(a)}(t),Z^{(a)}(t),t\geq
0\}$ satisfying the Skorokhod equation \bea\label{SEqcts} Z^{(a)}(t)
& = & H^{(a)}(t) + RY^{(a)}(t),~t\geq 0,\eea such that
$Y^{(a)}(0)=0,Z^{(a)}(0)=a,$ $Z^{(a)}(t)\in\overline{G},t\geq 0,$
each component of $Y^{(a)}(\c)$ is nondecreasing, and
$(Y^{(a)})_i(\c)$ can increase only when $(Z^{(a)})_i(\c)=0,$ in the
sense that \bea\label{Mincts} (Y^{(a)})_i(t)-(Y^{(a)})_i(s) & = &
\int_{(s,t]}I_{\{0\}}((Z^{(a)})_i(u))d(Y^{(a)})_i(u), \eea for
$0\leq s<t.$ Thanks to (H1), this Skorokhod problem has a unique
solution; the solution also has desired optimal property as in the
discrete time.

Put $\Delta Y^{(a)}(t)=Y^{(a)}(t)- Y^{(a)}(t-),~t\geq 0.$ We say
that ss-ruin occurs for the process $Z^{(a)}(\c)$ if $\Delta
Y^{(a)}(t) \gg 0$ for some $t>0;$ by \eqref{Mincts}, note that
ss-ruin means that every company has zero surplus and nonzero
marginal deficit at some time $t.$ Other two notions of ruin can be
similarly defined as in Subsection 3.2.  As $c\gg 0$ and
$A_{\ell}>0,$ note that each coordinate of $H^{(a)}(\c),$ and hence
that of $Z^{(a)}(\c)$ is strictly increasing between claim arrivals;
so  $(Y^{(a)})_i$ can change only at a claim arrival time, for any
$i.$ In particular, only at a claim arrival time ruin can occur.
Therefore, as in the one dimensional case, for studying ruin
problem, it is sufficient to consider these processes at claim
arrival times. For $n=1,2,\cs$ let $T_n=A_1+A_2+\cs+A_n$ denote the
arrival times, and let $T_0=0;$ set $H^{(a)}_n=H^{(a)}(T_n),n\geq
0.$ Then $H^{(a)}_n,n\geq 0$ is a random walk in $\R^d$ starting at
$a.$ Put
$Y^{(a)}_n=Y^{(a)}(T_n),Z^{(a)}_n=Z^{(a)}(T_n),~n=0,1,2,\cs$ Then
$\{Z^{(a)}_n\}$ is the associated regulated/ reflected random walk
with $\{Y^{(a)}_n\}$ as the corresponding pushing process. Note that
$\Delta Y^{(a)}_n=\Delta Y^{(a)}(T_n),n\geq 1.$ So, ss-ruin for the
process $\{Z^{(a)}(t),t>0\}$ occurs when and only when ss-ruin for
the regulated random walk $\{Z^{(a)}_n,n\geq 1\}$ occurs. Moreover,
the respective probabilities of ss-ruin in finite time also
coincide. Similar comments apply also to other notions of ruin.
Another consequence of the above discussion is \bea\label{ytyn} 
M &=& \lim_{n\r\iy}\{Y^{(0)}_n:\Delta Y^{(0)}_n \gg 0\} \nonumber \\
&=& \lim_{t\r\iy} \{Y^{(0)}(t):\Delta Y^{(0)}(t) \gg 0\} . \eea 
Assume now that (H1)-(H6),(H8) hold.
Therefore Theorem~\ref{Th:duality3} can be used to conclude that
\bea\label{renewruinprob} \P( \varrho^{(a)}_{ss} < \iy) &=&
\P(M \gg R^{-1}a ), \eea giving the probability of
ruin in finite time for the regulated/ reflected process
$\{Z^{(a)}(t):t\geq 0\}$ corresponding to initial capital vector $a,$
in terms of an asymptotic functional of the pushing process
$\{Y^{(0)}(t):t\geq 0\}$ corresponding to zero initial capital
exceeding threshold $R^{-1}a.$  }

\end{Example}

\begin{Example}\label{Ex:onedim} $1-$dimensional problem revisited:
{\rm Take $d=1.$ The classical renewal risk process is then given by
the one-dimensional analogue $H^{(a)}(\c)$ of \eqref{Hcts}. In this
case the scalar $R=1.$ Consider the one-dimensional Skorokhod
problem in the half line $[0,\iy)$ for $H^{(a)}(\c);$ let
$Z^{(a)}(\c),Y^{(a)}(\c)$ denote, respectively, the regulated/
reflected process and the pushing process. If the claim size
distribution is continuous, then it is easily seen that the events
$\{H^{(a)}(t)< 0$ for some $t>0\}$ and $\{\Delta Y^{(a)}(t)>0$ for
some $t>0\}$ coincide with probability $1.$ So the classical notion
of ruin probability and all 3 notions of ruin probability discussed
here are the same with probability $1.$ In particular, by 
Corollary~\ref{Cor:limdist} and 
\eqref{renewruinprob}, it follows that \bea\label{onedim}
\P(Y^{(0)}_{\iy} \gg R^{-1}a ) & = & \P(H^{(a)}(t) <0~{\rm for}~{\rm
some}~t>0). \eea For example, if claim size distribution is
Pareto$(\alpha)$ then $Y^{(0)}_{\iy}$ has Pareto$(\alpha-1)$
distribution. Thus the rich reservoir of results from ruin theory
(see \cite{AA2010,EKM1997,RSST1999} can be used to get information
on the asymptotics of the pushing process of Skorokhod problem.
}
\end{Example}

\begin{Example}\label{Ex:normalrefl} {\rm Suppose
$\P(X_{\ell} \gg x)>0,~\forall x\in\overline{G},~\ell =1,2,\cs$ Take
$R=I,$ the identity matrix. Then (H9) given in Remark~\ref{Rem:reH2}
is satisfied; clearly (H1) holds. If  (H3)-(H5),(H7),(H8) also hold,
then by Remark~\ref{Rem:reH2}, our analysis is applicable in this
case as well. This situation corresponds to the companies operating
without the risk diversifying treaty, claim arrival times being the
same for all companies, but claim sizes may be dependent; that is,
the $d$ companies jointly take care of different (possibly
dependent) components of the vector claims, or the companies take
care of each claim in nonzero proportions. Capital injection needed
by a company is completely taken care of by its shareholders. In
this set up, ruin of the network is the same as the $d-$dimensional renewal
risk process $H^{(a)}(\c)$ given by \eqref{Hcts} operating without
the risk diversifying treaty hitting the negative orthant. }
\end{Example}


\begin{thebibliography}{99}

\bibitem{AA2010} S.~Asmussen and H.~Albrecher: {\it Ruin Probabilities,}
(Second edition). World Scientific, Singapore, 2010.

\bibitem{APP2008} F.~Avram, Z.~Palmowski and M.~Pistorius: Exit problem of a
two-dimensional risk process from the quadrant: exact and asymptotic
results. {\it The Annals of Applied Probability} {\bf 18} (2008)
2421 -- 2449.

\bibitem{BG2008} N.~Bauerle and R.~Grubel: Multivariate risk processes with
interacting intensities. {\it Advances in Applied Probability} {\bf
40} (2008) 578 -- 601.

\bibitem{Bi2001} N.H.~Bingham: Random Walk and Fluctuation Theory.
In {\it Handbook of Statistics, Vol.19} (eds. D.N.~Shanbhag and C.R.~Rao),
Elsevier, Amsterdam, 2001.

\bibitem{BS1999} B.~Blaszczyszyn and K.~Sigman: Risk and duality in multidimensions. {\em Stoc. Proc. Appl. } {\bf 83} (1999) 331-356.

\bibitem{Bu1970} H.~Buhlman: {\it Mathematical Methods in Risk Theory.}
Springer-Verlag, Berlin-Heidelberg, 1970.

\bibitem{CYZ2003} W.-S.~Chan, H.~Yang and L.~Zhang: Some results on ruin probabilities in a two-dimensional risk model. {\it Ins. Math. Econ. } {\bf 32} (2003) 345-358

\bibitem{CM1991} H.~Chen and A.~Mandelbaum: Leontief systems, RBV's and
RBM's. In {\it Proceedings of Imperial College Workshop on Applied
Stochastic Processes,} (ed. M.H.A.~Davis and R.J.~Elliott), pp. 1-43.
Gordon and Breach, New York, 1991.

\bibitem{Co1998} J.~Collamore: First passage times of general
sequences of random vectors: a large deviations approach. {\it
Stochastic Processes and Applications} {\bf 78} (1998) 97-130.

\bibitem{CPS1992} R.W.~Cottle, J.S.~Pang and R.E.~Stone: {\it The
Linear Complementarity Problem.} Academic Press, New York, 1992.

\bibitem{DW2004} D.C.M.~Dickson and H.R.~Waters: Some optimal dividends problem.
{\it ASTIN Bulletin} 34 (2004) 49-74.

\bibitem{EKM1997} P.~Embrechts, C.~Kluppelberg and T.~Mikosch: {\it
Modelling Extremal Events for Insurance and Finance.} Springer,
Heidelberg, 1997.

\bibitem{F1969} W.~Feller: {\it An Introduction to Probability Theory and
its Applications, } Vol. II. Wiley-Eastern, New Delhi, 1969.

\bibitem{Ha1985} J.M.~Harrison: {\it Brownian Motion and Stochastic Flow
Systems.} Wiley, New York, 1985.

\bibitem{HR1981} J.M.~Harrison and M.I.~Reiman: Reflected Brownian
motion on an orthant. {\em Ann. Probab.} {\bf 9} (1981)
302-308.

\bibitem{K1996} O.~Kella: Stability and nonproduct form of stochastic fluid
networks with Levy inputs. {\em Ann. Appl. Prob. } {\bf 6} (1996) 186-199.

\bibitem{K1997} O.~Kella: Stochastic storage networks:
stationarity and the feedforward case. {\em J. Appl. Prob. } {\bf 34} (1997) 498-507.

\bibitem{KW1996} O.~Kella  and W.~Whitt: Stability and
structural properties of stochastic fluid networks. {\em J.
 Appl. Prob.} {\bf 33} (1996) 1169-1180.

\bibitem{Ma1989} A.~Mandelbaum: The dynamic complementarity problem.
{\em Lecture Notes,} Technion, Oct. 1989.

\bibitem{Pr1961} N.U.~Prabhu: On the ruin problem of collective risk theory.
{\em Ann. Math. Statis. } {\bf 32} (1961) 757-764.

\bibitem{R2006} K.~Ramanan: Reflected diffusions defined via
the extended Skorokhod map. {\em Elec. J. Probab. } {\bf 11}
(2006) 934-992.

\bibitem{Ra2000} S.~Ramasubramanian: A subsidy-surplus
model and the Skorokhod problem in an orthant. {\em Math.
 Oper. Res.} {\bf 25} (2000) 509-538.

\bibitem{Ra2006} S.~Ramasubramanian:  An insurance
network: Nash equilibrium. {\em Ins. Math.
Econ.} {\bf 38} (2006) 374-390.

\bibitem{Ra2011} S.~Ramasubramanian: Multidimensional
insurance model with risk-reducing treaty. {\em Stoc. Models}
{\bf 27} (2011) 363-387.

\bibitem{Ra2012} S.~Ramasubramanian: A multidimensional
ruin problem. {\em Comm. Stoc. Anal.} {\bf 6} (2012) 33-47.

\bibitem{Re1984} M.I.~Reiman: Open queueing networks in heavy traffic.
{\it Math. Oper. Res.} {\bf 9} (1984) 441-458.

\bibitem{RSST1999} T.~Rolski, H.~Schmidli, V.~Schmidt and
J.L.~Teugels: {\it Stochastic Processes for Insurance and Finance.}
Wiley, Chichester, 1999.

\bibitem{Se1981} E.~Seneta: {\it Non-negative matrices and Markov chains,}
(Second edition). Springer-Verlag, New York, 1981.

\bibitem{Sg1976} D.~Siegmund: The equivalence of absorbing and reflecting
barrier problems for stochastically monotone Markov processes.
{\em Ann. Probab. } {\bf 4} (1976) 914-924.

\bibitem{Sp1956} F.~Spitzer: A combinatorial lemma and its application to probability theory. {\em Trans. Amer. Math. Soc. } {\bf 82} (1956) 323-339.


\end{thebibliography}
\end{document}